\documentclass[11pt]{article}
\usepackage{amsmath,amsfonts,amsthm,amscd,amssymb,graphicx}
\numberwithin{equation}{section}

\newtheorem{theorem}{Theorem}[section]

\newcommand{\cO}{\mathcal{O}}

\def\beq{\begin{equation}}
\def\eeq{\end{equation}}
\def\bb1{{1\!\!1}}
\def\R{\mbox{Re }}
\def\I{\mbox{Im }}

\def\rit{{\Bbb R}}

\begin{document}


\title{Spectral stability of Prandtl boundary layers: 
\\ an overview}

\author{Emmanuel Grenier\footnotemark[1]
 \and Yan Guo\footnotemark[2] \and Toan T. Nguyen\footnotemark[3]
}

\date\today

\maketitle

\begin{abstract}
In this paper we show how the stability of Prandtl boundary layers is linked to the stability of shear flows
in the incompressible Navier Stokes equations. We then recall classical physical instability results, and give a short educationnal 
presentation of the construction of unstable modes for Orr Sommerfeld equations. We end the paper with a conjecture 
concerning the validity of Prandtl boundary layer asymptotic expansions. 
 \end{abstract}

\renewcommand{\thefootnote}{\fnsymbol{footnote}}

\footnotetext[1]{Equipe Projet Inria NUMED,
 INRIA Rh\^one Alpes, Unit\'e de Math\'ematiques Pures et Appliqu\'ees., 
 UMR 5669, CNRS et \'Ecole Normale Sup\'erieure de Lyon,
               46, all\'ee d'Italie, 69364 Lyon Cedex 07, France. Email: egrenier@umpa.ens-lyon.fr}

\footnotetext[2]{Division of Applied Mathematics, Brown University, 182 George street, Providence, RI 02912, USA. Email: Yan\underline{~}Guo@Brown.edu}

\footnotetext[3]{Department of Mathematics, Penn State University, State College, PA 16803. Email: nguyen@math.psu.edu.}


\section{Introduction}

This paper is motivated by the study of the inviscid limit of Navier Stokes equations in
a bounded domain. Let $\Omega$ be a subset of $\rit^2$ or $\rit^3$, and let us consider
the classical incompressible Navier Stokes equations in $\Omega$, posed on the velocity
field $u^\nu$,
\beq \label{NS1}
\partial_t u^\nu + \nabla ( u^\nu \otimes u^\nu) + \nabla p^\nu - \nu \Delta u^\nu = 0,
\eeq
\beq \label{NS2}
\nabla \cdot u^\nu = 0,
\eeq 
with no--slip boundary condition 
\beq \label{NS3}
u^\nu = 0 \qquad \hbox{ on }\quad \partial \Omega .
\eeq
As the viscosity $\nu$ goes to $0$, we would expect to recover incompressible Euler equations
\beq \label{Eu1}
\partial_t u^0 + \nabla (u^0 \otimes u^0) + \nabla p^0 = 0,
\eeq
\beq \label{Eu2}
\nabla \cdot u^0 = 0,
\eeq
with boundary condition
\beq \label{Eu3}
u^0 \cdot n = 0 \qquad \hbox{ on }\quad \partial \Omega,
\eeq
where $n$ is the outer normal to $\partial \Omega$. Throughout the paper, for the sake of presentation, we shall assume that $\Omega$ is the two-dimensional half space with $z\ge 0$.  

\medskip

The no-slip boundary condition (\ref{NS3}) is the most difficult condition to study the inviscid limit problem.
It is indeed the most classical one and the genuine one, historically considered
in this framework by the most prominent physicists including Lord Rayleigh, W. Orr, A. Sommerfeld,
W. Tollmien, H. Schlichting, C.C. Lin, P. G. Drazin, W. H. Reid, and L. D. Landau, among many others. See for example the physics books on the subject: Drazin and Reid \cite{Reid} and Schlichting \cite{Schlichting}. 
If the boundary condition (\ref{NS3}) is replaced by the Navier (slip) condition, boundary layers, though sharing the same thickness of $\sqrt \nu$, have much smaller amplitude (of an order $\sqrt \nu$, instead of order one of the Prandtl boundary layer), and are hence more stable (the smaller the boundary layer is, the more stable it is). We refer for instance to \cite{Iftimie, Sueur, Masmoudi}
for very interesting mathematical studies of boundary layers under the Navier boundary conditions.

\medskip

It is then natural to ask whether $u^\nu$ converges to $u^0$ as $\nu \to 0$ with the no-slip boundary condition \eqref{NS3}. This question
appears to be very difficult and widely open in Sobolev spaces, mainly because the boundary condition
changes between the Navier Stokes and Euler equations. Precisely, the tangential velocity
vanishes for the Navier Stokes equations, but not for Euler. In the limiting process
a boundary layer appears, in which the tangential velocity quickly goes from the Euler
value to $0$ (the value of the Navier-Stokes velocity on the boundary).

\medskip

The boundary layer theory was invented by Prandtl back in 1904 (when the first boundary layer equation
was ever found). Prandtl assumes that the velocity in the boundary layer depends on $t$, $x$ and on
a rescaled variable 
$$
Z = {z \over \lambda}
$$
where $\lambda$ is the size of the boundary layer. We therefore make the following Ansatz, within the 
boundary layer,
$$
u^\nu(t,x,z) = u^P(t,x,Z) + \lambda u^{P,1}(t,x,Z) + ... .
$$
Let the subscript $1$ and $2$ denote horizontal and vertical components of the velocity, respectively. 
The divergence free condition \eqref{NS2} then gives
$$
\Big(\partial_x u^P_1 + \lambda \partial_x u^{P,1}_1 + \cdots \Big)
+ \Big( \lambda^{-1} \partial_Z u^P_2 + \partial_Z u^{P,1}_2 + \cdots \Big) = 0,
$$
which by matching the respective order in the limit $\lambda\to 0$ in particular yields
\begin{equation}\label{Pr0}
\partial_Z u^P_2 = 0, \qquad 
\partial_x u^P_1 + \partial_Z u^{P,1}_2 = 0.
\end{equation}
As $u^\nu$ vanishes at $z=0$, this implies that 
$$
u_2^P = 0,
$$
identically: the vertical velocity in the boundary layer is of order $O(\lambda)$.


Now, the Navier-Stokes equation \eqref{NS1} on the horizontal speed gives 
\beq \label{Pr1}
\partial_t u^P_1 + u^P_1 \partial_x u^P_1 
+ u^{P,1}_2 \partial_Z u^P_1 -  \partial^2_{ZZ} u^P_1 = - \partial_x p,
\eeq
provided that we choose
$$
\lambda = \sqrt{\nu} .
$$
Next, to leading order, the equation on the vertical speed reduces to 
\begin{equation}\label{Pr2}
\partial_Z p = 0.
\end{equation}
Hence the leading pressure $p$ depends only on $t$ and $x$, and is given by the
pressure at infinity, namely by the pressure of Euler flow in the interior of the domain.  
As for boundary conditions, we are led to impose
\beq \label{Pr3}
u^P_1 = u^{P,1}_2 = 0 \qquad \hbox{ on }\quad  \{Z = 0\}
\eeq
and 
\beq \label{Pr4}
u_1^P(t,x,Z)  \to u^E_1(t,x,0) \qquad \hbox{ as }Ê\quad Z \to + \infty,
\eeq
where $u^E(t,x,0)$ denotes the value of the Euler flow in the interior of the domain (away from
the boundary layer). The set of equations \eqref{Pr0}-\eqref{Pr4} is called the Prandtl boundary layer equations. 

\medskip


A natural question then arises: {\em can we justify that $u^\nu$ is the sum of an Euler
part $u^E$ plus the Prandtl boundary layer  correction $u^P$ ?}

\medskip

The first problem is to prove existence of solutions for the Prandtl equation. This is difficult since
whereas $u^P_1$ satisfies the simple transport equation with the degenerate diffusion, 
$u^{P,1}_2$ satisfies no prognostic equation, and can only be recovered, using
$$
u^{P,1}_2(t,x,Z) = - \int_0^Z \partial_x u^P_1(t,x,Z) dZ.
$$
Hence $u^{P,1}_2$ is the vertical primitive of an horizontal derivative.
This leads to the loss of one derivative in the estimates. 
In the analytic framework, it is possible to 
control one loss of derivative: the Prandtl equation is well posed for small times; see \cite{SC1,KV}. See also \cite{GM} for the construction of Prandtl solutions in Gevrey classes. The existence of Prandtl solutions
in Sobolev spaces is delicate. Oleinik \cite{Oleinik} was the first to establish the existence of smooth solution in finite time
provided that the initial tangential velocity $u_1^P(0,x,Z)$ is monotonic in $Z$. Monotonicity plays a crucial
role in its proof and makes it possible for the existence via special transformations; see also recent works \cite{Alex, MW, KMVW} where the solution is constructed via delicate energy estimates. Then E and Engquist \cite{EE} proved that Prandtl layer may blow up in finite time.
More recently, G\'erard-Varet and Dormy \cite{DGV} proved that the Prandtl equation is linearly ill-posed in Sobolev
spaces if Oleinik's monotonicity assumption is violated. 


\medskip

Concerning the justification of boundary layers, the analytic framework has been investigated in full details by Sanmartino and Caflisch in \cite{SC1,SC2}.
They prove that, with analytic assumptions on the initial data, the Navier Stokes solution can be described asymptotically as the sum of an Euler solution in the interior and a Prandtl boundary layer correction. Recently, the author \cite{Mae} was able to prove the $L^\infty$ convergence under the assumption that the initial vorticity is away from the boundary.

\medskip

These results in particular prove that Prandtl boundary layers are the right expansion, since 
if there is an expansion, it should be true for analytic
functions, and therefore it must involve Prandtl layers. Therefore, we have no alternative asymptotic
expansions.

\medskip

However, analytic regularity is a very strong assumption. It mainly says that
there are no high frequencies in the fluid (energy spectrum of noise
decreases exponentially as the spatial frequency goes to infinity). In physical cases however there is always
some noise, which is not so regular (energy decreases like an inverse power
of the spatial frequency). 


\medskip

Let us from now on consider Sobolev regularity. 
In general, it does not appear to be possible to prove that Navier Stokes solutions behave
like Euler solutions plus a Prandtl boundary layer correction if we seek for global-in-time 
results or if initially the boundary layer profile has an inflection point or the profile is not monotonic; see \cite{Gr1,GN}. Though, it leaves open that whether this expansion is possible for small time and monotonic initial
profiles with no inflection point in the boundary layer. The aim of this program is to discuss this question
in the case of shear flows, where the limiting Euler equation is trivial $u^E(t,x,z) = U_\infty$
(constant flow). Of course, a non-convergent result in this particular case would indicate that the expansion is not possible in general.


\section{Inside the boundary layer}


As mentioned earlier, it is crucial to understand what happens inside the boundary layer, which is of the size $\sqrt \nu$. Prandtl chooses an anisotropic change of variables
$$
T = t, X = x, Z = {Z \over \sqrt{\nu}}.
$$
However, a natural tendency of fluids is to create vortices, and vortices tend to be isotropic (comparable
sizes in $x$ and $z$). Vortices also evolve within times of order of their size. Hence it is more natural to 
introduce an isotropic change of variables
$$
(T,X,Z) = {1 \over \sqrt{\nu}} (t,x,z) .
$$
In these new variables, the system of equations (\ref{NS1}), (\ref{NS2}), and (\ref{NS3}) turns to
\beq \label{NS1b}
\partial_T u^\nu + \nabla  ( u^\nu \otimes u^\nu) + \nabla p^\nu - \sqrt{\nu} \Delta u^\nu = 0,
\eeq
\beq \label{NS2b}
\nabla \cdot u^\nu = 0,
\eeq 
with no-slip boundary condition 
\beq \label{NS3b}
u^\nu = 0 \qquad \hbox{ on } \quad \partial \Omega .
\eeq
These equations are again the Navier Stokes equations where the viscosity $\nu$ has
been replaced by $\sqrt{\nu}$. These equations admit particular solutions of the form
\beq \label{shear}
u^\nu(T,X,Z) = U^P(\sqrt{\nu} T, Z) 
\eeq
with $$
U^P(t,Z) = (U^P_s(t,Z), 0),
$$
where $U^P_s$ satisfies the scalar heat equation
\beq \label{heat1}
\partial_t U^P_s(t,Z) - \partial_Z^2 U^P_s(t,Z) = 0,
\eeq 
with boundary condition
\beq \label{heat2}
U^P_s(t,0) = 0 .
\eeq
The particular solution $U^P$ is called the shear flow or shear profile. Note that $U^P(t,Z)$ is also a particular solution of the Prandtl equations, since for
shear layer profiles, the Prandtl equations and Navier Stokes equations simply 
reduce to the same heat equation. The existence of solutions to Prandtl equation is of course
trivial in this particular case.

\medskip

{\em However, do we still have convergence for small Sobolev perturbations of such profiles?}

\medskip

Namely, let us consider initial data of the form
\beq \label{init}
u^\nu(0,X,Z) = (U^P_s(0,Z), 0) + v^\nu(0,X,Z),
\eeq
where $v^\nu$ is initial perturbation that is small in Sobolev spaces. Do we still have convergence
of $u^\nu(T,X,Z)$ to $U^P(\sqrt{\nu} T, Z)$, for $T>0$ ?

On bounded time intervals $0 < T < T_0$ ($T_0$ is fixed and independent on $\nu$),
the convergence is true and can be seen easily through classical $L^2$ energy estimates. However we are interested
by results on time intervals of the form $0 < T < T_0 / \sqrt{\nu}$ (that is, a uniform time in the original variable $t = \sqrt \nu T$). On such a long
interval in the rescaled variables, the classical $L^2$ energy estimates are useless.  The problem is to know
whether small perturbations of the limiting Prandtl profile can grow in a large time. {\em This is  a stability problem
for a shear profle for Navier Stokes equations.}

\medskip

The first step is to look at the linearized stability of the shear layer $U^P(\sqrt \nu T,Z)$. Let us freeze the time dependence in this shear profile, and study the stability of the time-independent profile $U^P(0,Z)$. The linearized Navier Stokes equations near $U^P(0,Z)$ then read
\beq \label{NS1l}
\partial_T v^\nu +  ( U^P \cdot \nabla) v^\nu + (v^\nu \cdot \nabla) U^P 
+ \nabla q^\nu - \sqrt{\nu} \Delta v^\nu = 0,
\eeq
\beq \label{NS2l}
\nabla \cdot v^\nu = 0,
\eeq 
with no-slip boundary condition 
\beq \label{NS3b}
v^\nu = 0 \qquad \hbox{ on }\quad  \partial \Omega .
\eeq

If all the eigenvalues of this spectral problem have non-positive real parts, then it is likely
that $v^\nu$ remains bounded for all time, and that this is also true for the linearization near
the time-dependent profile $U^P(\sqrt \nu T,Z)$ and also true for the nonlinear Navier Stokes equations. In this case, we could
expect convergence from Navier Stokes to Euler with a Prandtl correction.

If one eigenvalue has a positive real part, then there exists a growing mode of the form
$V e^{c(\nu) T}$, with $\R c(\nu) > 0$. The time scale of instability $1 / \R c(\nu)$ must
then be compared with $1 / \sqrt{\nu}$. If $\R c(\nu) \ll \sqrt{\nu}$, then instability appears
in very large time, much larger than $T_0 / \sqrt{\nu}$ and convergence may hold.
On the contrary if $\R c(\nu) \gg \sqrt{\nu}$, then instability is strong and occurs much before
$T_0 / \sqrt{\nu}$. In this latter case, it is then likely that such an instability occurs
for $U^P(0,Z)$ and that it might not possible to prove convergence of Navier Stokes to Euler plus a Prandtl layer
in supremum norm or strong Sobolev norms.

The study of Prandtl boundary layer is therefore closely linked to the question of the spectral 
stability of shear profiles for Navier Stokes equation with $\sqrt{\nu}$ viscosity, 
and more precisely to the comparison of $\Re(c)$ with respect to $\sqrt{\nu}$.


\section{Spectral problem}



\subsection{Orr Sommerfeld and Rayleigh equations}


The analysis of the spectral problem is a very classical issue in fluid mechanics.
A huge literature is devoted to its detailed study. We in particular refer to
\cite{Reid, Schlichting} for the major works of Tollmien, C.C. Lin, and Schlichting.
The studies began around 1930, motivated by the study of the boundary layer
around wings. In airplanes design, it is crucial to study the boundary layer
around the wing, and more precisely the transition between the laminar and turbulent
regimes, and even more crucial to predict the point where boundary layer
splits from the boundary. A large number of papers has been devoted to 
the estimation of the critical Rayleigh number of classical shear flows 
(Blasius profile, exponential suction/blowing profile, etc...). 

Let us go further in detail in the case of two dimensional spaces. The first step is to make a Fourier transform
with respect to the horizontal variable, and a Fourier transform with respect to time variable 
on the to  stream function $\phi$. This leads to the following form for perturbations
$v^\nu$
\begin{equation}\label{def-stream}v^\nu = \nabla^\perp \psi = (\partial_Z, -\partial_X)\psi 
,\qquad \psi(T,X,Z) := \phi (Z) e^{i \alpha (X - cT) }.
\end{equation}
Putting this Ansatz in (\ref{NS1l}), we get the classical Orr--Sommerfeld equation
\beq \label{OS1}
{1 \over i \alpha R} (\partial_Z^2 - \alpha^2)^2 \phi 
= (U-c) ( \partial_Z^2 - \alpha^2) \phi  - U'' \phi  
\eeq
with boundary conditions
\beq \label{OS2}
\alpha \phi  = \partial_Z \phi  = 0 \qquad\hbox{ at } Z = 0
\eeq
and
\beq \label{OS3}
 \phi  \to 0\qquad  \hbox{ as } Z \to + \infty.
\eeq
Here $R = 1 / \sqrt{\nu}$ is the Reynolds number (to our rescaled equations) and $U = U_s^P(0,Z)$
is the shear profile  introduced in \eqref{heat1} and \eqref{heat2}.
The spectrum of (\ref{OS1}) clearly  depends on $\alpha$ and $R$.

\medskip

As $R \to \infty$, or rather $\alpha R \to \infty$, the Orr--Sommerfeld equations formally reduce to the so-called
Rayleigh equation
\beq \label{Rayleigh1}
 (U-c) ( \partial_Z^2 - \alpha^2) \phi  = U'' \phi  
\eeq
with boundary conditions
\beq \label{Rayleigh2}
 \phi  = 0\qquad  \hbox{ at } Z = 0
\eeq
and
\beq \label{Ra3}
\phi  \to 0\qquad  \hbox{ as } Z \to + \infty.
\eeq
The Rayleigh equation describes the stability of the shear profile $U$ for  Euler equations.
The spectrum of Orr Sommerfeld is a perturbation of the spectrum of Rayleigh equation. It is therefore natural to first study the Rayleigh equation.

\medskip

Stability of the Rayleigh problem depends on the profile. For some profiles, all the eigenvalues
are imaginary, and for some others there exist unstable modes. There are various
criteria to know whether a profile is stable or not, including classical Rayleigh inflection point
and Fjortoft criteria. We shall recall these two criteria in the next subsection.


\subsection{Classical stability criteria}
The first criterium is due to Rayleigh. 
\medskip

{\bf Rayleigh's inflexion-point criterium (Rayleigh \cite{Ray}).}  {\em A necessary condition for instability is that the basic
velocity profile must have an inflection point.}

\medskip

The criterium can easily be seen by multiplying by $\bar \phi/(U-c)$ to the Rayleigh equation \eqref{Rayleigh1} and using integration by parts. This leads to 
\begin{equation}\label{id-Ray}\int_0^\infty (|\partial_Z \phi|^2 + \alpha^2 |\phi|^2) \;dZ+ \int_0^\infty \frac{U''}{U-c}|\phi|^2 \;dZ=0,\end{equation}
whose imaginary part reads 
\begin{equation}\label{id-Ray1}\I c \int_0^\infty \frac{U''}{|U-c|^2}|\phi|^2 \;dZ =0.\end{equation}
Thus, the condition $\I c>0$ must imply that $U''$ changes its sign. This gives the Rayleigh criterium. 

\medskip

A refined version of this criterium was later obtained by Fjortoft (1950) who proved 
\medskip

{\bf Fjortoft criterium \cite{Reid}.} {\em A necessary condition for instability is that $U'' (U - U(z_c))<0$ somewhere in the flow, where $z_c$ is a point at which $U''(z_c) =0$.}

\medskip
To prove the criterium, consider the real part of the identity \eqref{id-Ray}: 
$$\int_0^\infty (|\partial_Z \phi|^2 + \alpha^2 |\phi|^2) \;dZ+ \int_0^\infty \frac{U'' (U-\R c)}{|U-c|^2}|\phi|^2 \;dZ=0.$$
Adding to this the identity 
$$ (\R c - U(z_c))\int_0^\infty \frac{U''}{|U-c|^2}|\phi|^2 \;dZ =0,$$
which is from \eqref{id-Ray1}, we obtain 
$$ \int_0^\infty \frac{U'' (U-U(z_c))}{|U-c|^2}|\phi|^2 \;dZ=-\int_0^\infty (|\partial_Z\phi|^2 + \alpha^2 |\phi|^2) \;dZ <0,$$
 from which the Fjortoft criterium follows.

\medskip

\subsection{Unstable profiles for Rayleigh equation}

If the profile is unstable for the Rayleigh equation, then there exist $\alpha$ and an eigenvalue $c_\infty$ with
$\I c_\infty > 0$, with corresponding eigenvalue $\phi _\infty$. 
We can then make a perturbative analysis to construct an eigenmode $\phi _R$ of
the Orr-Sommerfeld equation with an eigenvalue $\I c_R >0$ for any large enough $R$.

The main point is that $\partial_Z \phi _R$ vanishes on the boundary whereas $\partial_Z\phi _\infty$
does not necessarily vanishes. We therefore need to add a boundary layer to
correct $\phi _\infty$. This boundary layer comes from the balance between the terms
$\partial_Z^4 \phi  / \alpha R$ and $U \partial_Z^2 \phi $ of (\ref{OS1}) and is therefore
of size
$$
\sqrt{ U_0 \over \alpha R} = \cO( R^{-1/2}) = \cO(\nu^{1/4}).
$$
In original $t,x,y$ variables, this leads to a boundary layer of size $O(\nu^{3/4})$.
In the limit $\nu \to 0$, two layers appear: the Prandtl layer of size $\sqrt{\nu}$
and a so-called viscous sublayer of size $\nu^{3/4}$.
This sublayer has an exponential profile in $Z/ \nu^{1/4}$.
The existence and study of the viscous sublayer is a classical issue in physical fluid mechanics.

When $\phi _R$ is constructed and corrected by this sublayer, it in fact still does not satisfy (\ref{OS1}), but it does satisfy the Orr-Sommerfeld boundary conditions exactly and the Orr-Sommerfeld equation
up to an error with size of $\cO(1/R)$. By perturbative arguments we can prove
\beq \label{perturb1}
c_R = c_\infty + \cO(R^{-1}) .
\eeq
Next, starting from $\phi _R$, we can then construct unstable modes for the linearized Navier
Stokes equations, and even get instability results in strong norms for the nonlinear Navier
Stokes equations. This has been carried out in detail by E. Grenier in \cite{Gr1}.


\subsection{Stable profiles for Rayleigh equation}


Some profiles are stable for the Rayleigh equation; in particular, shear profiles
without inflection points from the Rayleigh's inflexion-point criterium. For stable profiles, all the spectrum of the Rayleigh equation
is imbedded on the axis: $\I c = 0$.
At a first glance, we may believe that (\ref{perturb1}) still holds true, which would mean
that any eigenvalue of the Orr-Sommerfeld would have an imaginary part $\I c_R$ of
order $\cO(R^{-1}) = \cO(1 / \sqrt{\nu})$. This would mean that perturbations would
increase slowly, and only get multiplied by a constant factor for times $t$ of order
$T_0 / \sqrt{\nu}$. In this case we might hope to obtain the convergence from Navier-Stokes
to Euler and Prandtl equations. {\em However, this is not the case, and $\I c_R$ appears to be much larger.} Let us detail 
now this point.

\medskip

The main point is that in the case of a stable profile, there exists an eigenmode $\phi_\infty$ with corresponding
eigenvalue $c_\infty$ which is {\it small} 	and {\it real}. Therefore there exists some $z_c$ such that
$$
U(z_c) = c_\infty .
$$
Such a $z_c$ is called a critical layer. As $z_c$, $U(z) - c_\infty$ vanishes, hence Rayleigh
equation is singular
\beq \label{RayR}
(\partial_Z^2 - \alpha^2) \phi = {U'' \over U - c_\infty} \phi .
\eeq
Therefore when $R$ goes to infinity, for $z$ near $z_c$, Orr Sommerfeld degenerates from a fourth
order elliptic equation to a singular second order equation. At $z = z_c$, all the derivatives disappear
as $R$ goes to infinity, and we go from a fourth order equation to a "zero order" one. The limit is 
therefore very singular, and as a matter of fact $\Im c(R)$ is much larger than expected.

\medskip

Let us go on with the analysis of Rayleigh equation.
The Rayleigh equation (without taking care of boundary 
conditions) admits two independent solutions $\phi _1$ and $\phi _2$, one smooth
$\phi _1$ which vanishes at $z_c$ and another $\phi _2$ which is less regular near
$z_c$. Using (\ref{RayR}) we see that $\phi_2''$ behaves like $O(1 / Z - z_c)$
near $z_c$. Hence $\phi_2$ behaves like $(Z-z_c) \log(Z-z_c)$ near  $z_c$.
Therefore, the eigenvector $\phi _\infty$ is of the form
\beq \label{form}
\phi _\infty = P_1(Z) + (Z - z_c) \log(Z-z_c) P_2(Z)
\eeq
where $P_1$ and $P_2$ are smooth functions, with $P_1(z_c) = 0$.

\bigskip

If we try to make a perturbation analysis to get $\phi _R$ out of $\phi _\infty$, we then face two
difficulties. First, we have to correct $\phi _\infty$ in order to satisfy $\phi _\infty = 0$ at $Z = 0$. 
But there is another much more delicate difficulty.
As $\phi _\infty$ is not smooth at $Z = z_c$ it is not a good approximation of $\phi _R$ near
$z_c$. In particular $(\partial_Z^2 - \alpha^2)^2 \phi _\infty$ is too singular at $z_c$, of order
$O(1 / (Z - z_c)^3)$.

To find a better approximation, one notes that near the singular point $z_c$, the term $\partial_Z^4 \phi _R$ 
can no longer be neglected.
In fact, near this point $z_c$, $\partial_Z^4\phi  / i \alpha R$ must balance with $U'(z_c) (Z - z_c) \partial_Z^2 \phi $.
This leads  to the introduction of another boundary layer of size $(\alpha R)^{1/3}$, near $z_c$
 satisfying the equation
\beq \label{Airy}
\partial_Y^2 \Phi_R = Y \Phi_R, \qquad \Phi_R := \partial_Z^2 \phi _R
\eeq
where
$$
Y := {Z - z_c \over (\alpha R U'(z_c))^{1/3}}. 
$$
This layer is called {\it critical layer}.
Note that (\ref{Airy}) is simply the classical Airy equation. If we try to construct $\phi _R$ starting from 
$\phi _\infty$, we therefore have to involve Airy functions to describe what happens
near the critical layer. As a consequence, $(\alpha R)^{1/3}$ is an important parameter,
and similarly to the unstable case, we could prove
\beq \label{exp2}
c_R = c_\infty + \cO( (\alpha R)^{-1/3}) + \cO(R^{-1}) .
\eeq
Hence, the situation is very delicate. It has been intensively studied in the period
$1940 - 1960$ by many physicists, including Heisenberg, C.C. Lin, Tollmien, Schlichting, among others.
Their main objective was to compute the critical Reynolds number of shear layer flows,
namely the Reynolds number $R_c$ such that for $R > R_c$ there exists an unstable
growing mode for the Orr-Sommerfeld equation. Their analysis requires a careful study
of the critical layer.

\medskip

From their analysis, it turns out that there exists some $R_c$ (depening on the profile)
such that for $R > R_c$ there are solutions $\alpha(R)$, $c(R)$ and $\phi _R$ to the Orr-Sommerfeld equations with
$\I c(R) > 0$. Their formal analysis has been compared with modern numerical experiments
and also with experiments, with very good agreement. Note that physicists
are interested in the computation of the critical Reynolds number, since
any shear flow is unstable if the Reynolds number is larger than this critical
Reynolds number. {\em In this program, we are interested in the high Reynolds limit,
which is a different question.}
This limit is not a physical one, since any flow has a finite Reynolds
number, and not in any physical case can we let the Reynolds go to very very
high values. Physical Reynolds numbers may be large (of several millions
or billions), much larger than the critical Reynolds number, but despite their
large values, they are too small to enter the mathematical limit $R \to + \infty$
we are considering. Fluids would enter the mathematical asymptotic regime
if $R^{-1/7}$ or $R^{-1/11}$ (see below) are large numbers, which 
leads Reynolds numbers to be of order of billions of billions, much larger than
any physical Reynolds number! 

\medskip

It is thus important to keep in mind that the mathematical limit is not physically pertinent.
Physically, the most important phenomena are: the existence of a critical Reynolds
number (above which the shear flow is unstable), the transition from laminar
to turbulent boundary layers, the separation of the boundary layer from the boundary.
All of these occur near the critical Reynolds number, which is large, but not in the
asymptotic regime which we will now consider.

\medskip

The problem is now to study rigorously the asymptotic behavior of $\alpha$ and $c$ as
$R \to \infty$.
Let us present now some classical physical results. These results can be found, for example, in the book of Drazin and Reid \cite{Reid} or of H. Schlichting \cite{Schlichting}.

\medskip

For $R$ large enough there exists an interval $[\alpha_1(R),\alpha_2(R)]$
such that for every $\alpha$ in this interval there exists an unstable mode
with $\I c(R) > 0$. The asymptotic behavior of $\alpha_1$ and $\alpha_2$ depends
on the shear profile.

\medskip

\begin{itemize}

\item For plane Poiseuille flow (not a boundary layer): $U(z) = z^2 - 1$ for $0 < z < 1$.
In this case
$$
\alpha_1(R) \sim C_1 R^{-1 / 7}, \quad
\alpha_2(R) \sim C_2 R^{-1/11} .
$$
\item For boundary layer profiles:
$$
\alpha_1(R) \sim C_1 R^{-1 / 4}, \quad
\alpha_2(R) \sim C_2 R^{-1/6} .
$$
\item For the Blasius (a particular boundary layer) profile:
$$
\alpha_1(R) \sim C_1 R^{-1 / 4}, \quad
\alpha_2(R) \sim C_2 R^{-1/10} .
$$

\end{itemize}

\begin{figure}[t]
\centering
\includegraphics[scale=.4]{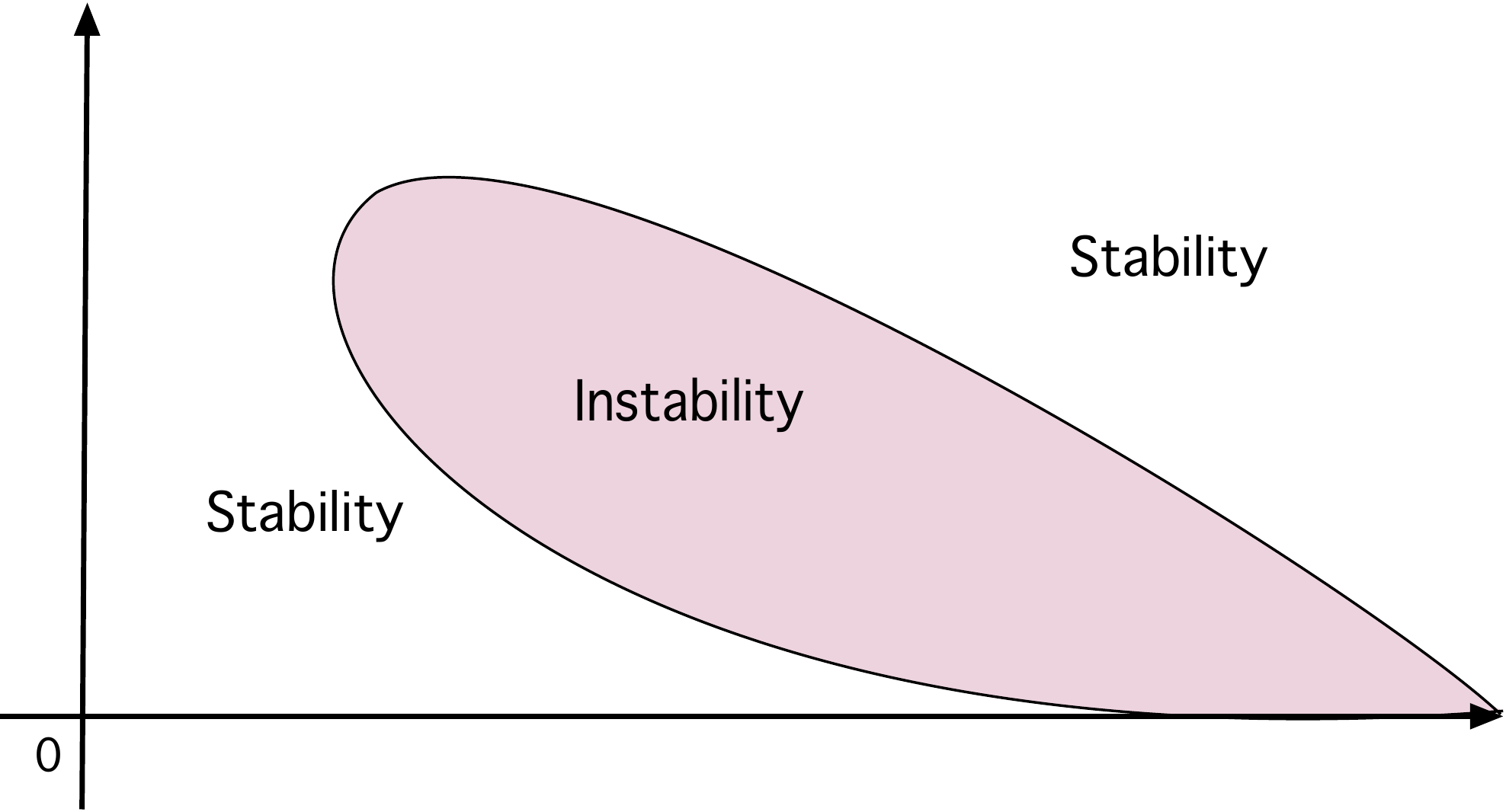}
\put(-20,1){$R^{1/5}$}
\put(-210,117){$\alpha^2$}
\put(-213,30){$\alpha_\mathrm{low}\approx R^{-1/4}$}
\put(-75,70){$\alpha_\mathrm{up}\approx R^{-1/10}$}
\caption{\em Illustrated are the marginal stability curves; see also \cite[Figure 5.5]{Reid}.}
\label{fig-shear}
\end{figure}
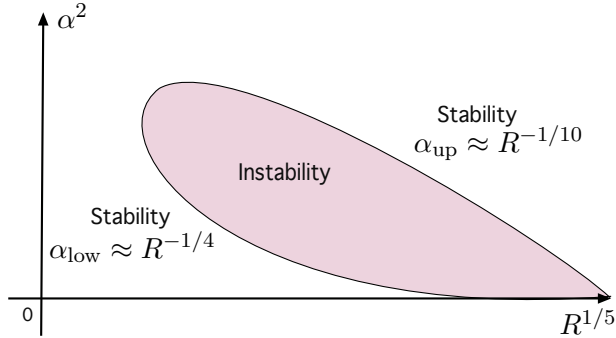

More precisely, in the $\alpha, R$ plane, the area where unstable modes exist is shown on figure \ref{fig-shear}. For small $R$, all the $\alpha$ are stable. Above some critical Reynolds number, there 
is a range $[\alpha_1(R),\alpha_2(R)]$ where instabilities occur. This instability area is bounded
by so called lower and upper marginal stability curves. 

Associated with the range of $\alpha(R)$, we have to determine the behavior of the eigenvalues
 $c(R,\alpha)$, or more precisely the imaginary part of $c(R,\alpha)$. 
 The complete mathematical justification of the construction of unstable modes will be detailled in companion
 papers. Here we juste want to present a quick and as simple as possible  construction of the unstable Orr 
 Sommerfeld modes. We will skip all difficulties and only focus on the backbone of the instability.

 
 \subsection{A sketch of the construction of unstable modes: the lower branch $\beta = 1/4$}
 
 
 We recall that there is no mathematically rigorous arguments of this section. We thus show the main
 ingredients of the instability, keeping under silence any other term. We assume that our profile
 $U(Z)$ is stable for Rayleigh equation. We focus on the lower
 marginal stability curve. In this case
 $$
 \alpha \sim A R^{-1/4} 
 $$
 and 
 $$
 \delta \sim (\alpha R)^{-1/3} \sim A^{-1/3} R^{-1/4} .
 $$
 Let us assume that $R$ is very large, and $\alpha$ very small. For small $\alpha$, Rayleigh
 equation is very close to
 \beq \label{RayT}
  (U - c) \partial_Z^2 \phi = U'' \phi
  \eeq
  which has an obvious solution
  $$
  \phi_1 = U - c.
  $$
 There exists another independent particular solution to (\ref{RayT}),  but it turns out that this second
 solution grows linearly as $Z$ increases, and may therefore be discarded. Note that
 $\phi_1$ is a smooth function and an approximate solution of Orr Sommerfeld.
 
 We next focus on Airy equation (\ref{Airy}). It has two particular fast decaying / growing solutions $\Phi_R = Ai$ and $Bi$.
 Only $Ai$ goes to $0$ as $Y$ goes to infinity, hence $Bi$ may be discarded.
 Let us denote by $Ai(1,Y)$ a primitive of $Ai$ and $Ai(2,Z)$ a primitive of $Ai(1,Y)$.
 Then
 $$
 \phi_3 = Ai(2,Y) = Ai(2, \delta^{-1} (Z - z_c) ) 
 $$
 is a particular solution of Airy, and an approximate solution of Orr Sommerfeld.
  Now we look for an eigenmode of Orr Sommerfeld which is a combination of $\phi_1$ and
 $\phi_3$ of the form:
 $$
 \phi = A \phi_1 + B \phi_3 .
 $$
  It has the good behavior as $Z$ goes to infinity. It remains to know whether we can find $A$ and $B$
  such that $\phi(0) = \partial_Z \phi(0) = 0$. This happens if the dispersion relation
  $$
  \phi_1(0) \partial_Z \phi_3(- \delta^{-1} z_c) = \partial_Z \phi_1(0) \phi_3(- \delta^{-1} z_c)
  $$
holds,  or equivalently if 
  \beq \label{disper}
  {\phi_1(0) \over \partial_Z \phi_1(0) } 
  =
  \delta {Ai(2,- \delta^{-1} z_c) \over Ai(1,- \delta^{-1} z_c)} .
  \eeq
  The left hand side of (\ref{disper}) is 
  $$
    {\phi_1(0) \over \partial_Z \phi_1(0) }  = {U_0 - c \over U_0'}
    $$
   and its imaginary part is simply $- \Im c / U_0'$, with $U_0' >0$. Here, $U_0 = U(0)$ and $U'_0 = U'(0)$.
   Now $Ai(2,Z) / Ai(1,Z)$ is the classical Tietjens function 
$$
T(Y)= {Ai(2,Y) \over Ai(1,Y) }.
$$
The main point is that $\Im T(Y)$ changes sign as $Y$ goes to infinity. It is positive for small $Y$ and
negative for large $Y$. As a consequence $\Im c$ changes sign as $z_c / \delta$ increases.
{\em This change of sign leads to the existence of unstable modes.}

It remains to link $z_c / \delta$ with $R$ and to prove that for $z_c / \delta$ goes to infinity as
$R$ increases. For this we have first to refine $\phi_1$. Namely $\phi_1$ does not go to $0$ as
$Z$ goes to infinity. For $\alpha > 0$, we may construct a solution $\phi_{1,\alpha}$ of the Rayleigh
equation which is a perturbation of $\phi_1$ and decreases like $\exp(-\alpha Z)$.
A classical perturbative analysis leads to
$$
\phi_{1,\alpha}(0) = U - c + \alpha { (U(\infty) - U_0)^2 \over U_0'} + ...
$$
Moreover as $Y$ goes to infinity,
$$
T(Y) \sim C Y^{-1/2} .
$$
Hence the dispersion relation takes the form
\beq \label{disper2}
{U_0 - c \over U_0'} + \alpha { (U(\infty) - U_0)^2 \over {U_0'}^2} + ...
\sim \delta (1 + | z_c  / \delta | )^{-1/2} .
\eeq
Assuming that $\alpha$ is much larger than the right hand side, which is the case if $A$ is
large enough, this gives that
$| U_0 - c |$  is of order $\alpha$, and hence $z_c$, defined by $U(z_c) = c$ is of order $\alpha$. Hence
$$
z_c / \delta \sim A^{4/3} .
$$
Therefore provided as $A$ increases, $\Im c$ changes from negative to positive values: there exists 
an threshold $A_{1c}$ such that $\Im c > 0$ is $A > A_{1c}$.
This ends our overview of the lower marginal curve.

 
 \subsection{A sketch of the construction of unstable modes: the upper branch $\beta = 1/6$}
 

The upper branch of marginal stability is more delicate to handle.
Roughly speaking, when the expansion of $\phi_{1,\alpha}$ involves $\phi_2$, independent
solution of Rayleigh equation which is singular like $(z - z_c) \log(z - z_c)$. This singularity
is smoothed out by Orr Sommerfeld in the critical layer. This smoothing involves second primitives
of solutions of Airy equation. As we take second primitives, a linear growth is observed 
(linear functions $\phi_R$ are obvious solution of (\ref{Airy})).  This linear growth gives an extra term
in the dispersion relation which can not be neglected when $\alpha \sim R^{-1/6}$. It has a stabilizing
effect and is responsible of the upper branch for marginal stability.
 

\section{Program}


The situation is well-known, physically speaking. However, to the best
of our knowledge, the formal analysis has never been justified mathematically.
Our ultimate goal is to prove the following conjecture:

\medskip

{\it
{\bf Conjecture:} generically, shear flows for Navier Stokes are linearly  unstable, 
and the Prandtl expansion is not valid in Sobolev spaces.}

\medskip

Let us now lay out our program to tackle the conjecture. 

\medskip

1. The first step is to construct unstable modes for the Orr-Sommerfeld equations
 as $R \to \infty$. This  requires a careful analysis of this singular
perturbation and a careful study of the behavior of the eigenvalues $c_R$. This leads to the proof
than generic shear layers are spectrally unstable.
%
%
%

More precisely, we will construct growing modes (those with $\I c >0$) for \eqref{OS1}-\eqref{OS3}
 when $R$ is large and $\alpha$ belongs to the interval $(\alpha_1(R),\alpha_2(R))$,
  with
   \begin{equation}\label{ranges-alpha}\alpha_1(R) = A_{1c}R^{-1/4}\qquad
     \mbox{and}\qquad  \alpha_2(R)= A_{2c} R^{-1/6}\end{equation}
    for some fixed constants $A_{1c}, A_{2c}$.  
    The curves $\alpha_j(R)$ are called lower and upper branches of the marginal (in)stability for the boundary layer $U$. 
That is, there is a critical constant $A_{1c}$ so that with $\alpha_1(R) = A_1 R^{-1/4}$, the imaginary part of $c$ turns from negative (stability) to positive (instability) when the parameter $A_1$ increases across $A_1=A_{1c}$. Similarly, there exists an $A_{2c}$ so that with $\alpha = A_2 R^{-1/6}$,
 $\I c$ turns from positive to negative as $A_2$ increases across $A_2 = A_{2c}$. In particular, we obtain instability of the profile in the intermediate zone: $\alpha \sim R^{-\beta}$ for $1/6<\beta<1/4$.

Our main result is as follows.
\begin{theorem}[Spectral instability of generic shear flows \cite{GGN1}]\label{theo-unstablemodes}~\\
Let $U(z)$ be a shear profile with $U'(0) \ne 0$ and satisfy 
$$
\sup_{z \ge 0} | \partial^k_z (U(z) - U_+) e^{\eta z} | < + \infty, \qquad k=0,\cdots ,4,
$$ for some constants $U_+$ and $\eta > 0$. 
There exists two constants $A_1$ and $A_2$ such that 
for $R$ large enough and for $\alpha = A R^{-\beta}$ with arbitrary $A$ if $1/6 < \beta < 1/4$
 or $A > A_1$ if $\beta = 1/4$ or $A < A_2$ if $\beta = 1/6$,
 there exist
$c(R)$  and  $\phi_R$ such that $\phi_R$ is an eigenfunction of the problem 
(\ref{OS1}), (\ref{OS2}), and (\ref{OS3}) with corresponding eigenvalue $c(R)$.  \\
More precisely, $\phi_R$ satisfies the boundary conditions (\ref{OS2})-(\ref{OS3})
and satisfies (\ref{OS1}). In addition, there holds the estimate 
$$
c(R) \quad \sim\quad  c_0 R^{\beta -1/2},
$$
for some constant $c_0$ independent on $R$, with $\I c_0 > 0 .$ In particular, the growth rate for the unstable modes is 
$$ \alpha \I c(R) \quad \sim\quad  R^{-1/2} .$$  
\end{theorem}

\subsubsection*{Remarks}

\begin{itemize}

\item[i)] The assumption $U'(0) \ne 0$ is technical. A similar analysis
could be fulfilled to allow the case $U'(0)=0$, with different (presumedly, more complicated) asymptotic behavior in the expansions.

\item[ii)]

The asymptotic behavior of the growth rate $\alpha \I c(R) \sim C_R R^{-1/2}$ holds in the {\it rescaled}
variables. In the original ones, this means that the unstable mode increases like
 $\exp( C t R^{1/2}) = \exp (C t /\nu^{1/4})$.  
 As a consequence,
one cannot expect stability in Sobolev norms for small perturbations of such shear flows.
Small perturbations will quickly increase in the time variable $t$ and may become of order
$1$ in a vanishing time (i.e., in a time that tends to zero as $\nu \to 0$). Therefore it is likely that slightly initially perturbed solutions of
Navier Stokes equations do not converge to the Prandtl equations as $\nu \to 0$.

\item[iii)]

It is worth noting that if we assume that the initial perturbation is analytic, then Fourier modes $\alpha / \sqrt{\nu}$
(in $x$ variables) are initially as small as  $\exp( - C / \sqrt{\nu})$.
Hence even if they grow fast, like 
$\exp(C_1 t / \nu^{1/4})$, 
they remain negligible as
long as $t < C  / (C_1\nu^{1/4})$. Therefore for small times, analytic perturbations remain negligible 
and we have convergence from Navier Stokes equation to Euler plus Prandtl for such
initial analytical data.

\end{itemize}

2. The second step is to prove linear instability. For a fixed viscosity, nonlinear instability follows from the spectral instability; see \cite{FPS} for arbitrary spectrally unstable steady states. However, in the vanishing viscosity limit, linear to nonlinear instability is a very delicate issue, primarily due to the fact that there are no available, comparable bounds on the linearized solution operator as compared to the maximal growing mode. Available analyses (for instance, \cite{Fri, Gr1}) do not appear applicable in the inviscid limit. In addition, boundary layers are shear layer profiles, which are time-dependent and are solutions of
the linear heat equation. In this case, even the proof of linear instability is no longer straightforward since the equation of the perturbation
changes with time. 

To get such a nonlinear instability result, we have to bound the resolvent of linearized Navier Stokes equations with fixed stationary 
profiles, and then treat the time-dependent profiles as small perturbations within a vanishing time in the inviscid limit. Getting bounds on the resolvent is
however highly technical, and we plan to follow the ideas developed by K. Zumbrun and coauthors; \cite{Zu1}. This problem 
will be investigated in a further work.

\medskip

Note that a similar analysis may be done for channel flows, including the classical plane Poiseuille flows. More precisely, we establish the following:

\begin{theorem}[Spectral instability of generic shear flows \cite{GGN2}]\label{theo-unstablemodes2}~\\
Let $U(z)$ be an arbitrary shear profile that is analytic and symmetric about $z=1$ with $U'(0) > 0$ and $U'(1) =0$. There exist 
 $\alpha_\mathrm{low}(R)$ and $\alpha_\mathrm{up}(R)$, there exists a critical Reynolds number $R_c$ so that for all $R\ge R_c$ and all $\alpha \in (\alpha_\mathrm{low}(R), \alpha_\mathrm{up}(R))$, there exist
a triple $c(R), \hat v(z; R), \hat p(z;R)$, with $\mathrm{Im} ~c(R) >0$, such that 
$$
v_R: = e^{i\alpha(y-ct) }\hat v(z;R), \quad
p_R: = e^{i\alpha(y-ct)} \hat p(z;R)
$$
solve the problem 
\eqref{NS1}-\eqref{NS2} with the no-slip boundary conditions. In the case of instability, there holds the following estimate for the growth rate of the unstable solutions:
$$ \alpha \I c(R) \quad \approx\quad  (\alpha R)^{-1/2},$$
as $R \to \infty$. In addition, the horizontal component of the unstable velocity $v_R$ is odd in $z$, whereas the vertical component is even in $z$. 
\end{theorem}


\end{document}